\documentclass[pdflatex,sn-mathphys-num]{sn-jnl}
\usepackage{amsmath,amssymb,amsfonts,amsthm}%
\usepackage{color,bbm}
\allowdisplaybreaks

\theoremstyle{thmstyleone}%
\newtheorem{theorem}{Theorem}[section]
\newtheorem{lemma}[theorem]{Lemma}%
\newtheorem{corollary}[theorem]{Corollary}%
\theoremstyle{thmstyletwo}%
\newtheorem{remark}[theorem]{Remark}%

\theoremstyle{thmstylethree}%

\begin{document}
\title{Maximal Noncompactness 

of Wiener-Hopf Operators}

\author*[1]{\fnm{Oleksiy} \sur{Karlovych}}
\email{oyk@fct.unl.pt}
\author[2]{\fnm{Eugene} \sur{Shargorodsky}}
\email{eugene.shargorodsky@kcl.ac.uk}
\equalcont{{These authors contributed equally to this work.}

\vspace{3mm}

To Sergei Grudsky on the occasion of his 70th birthday}
\affil*[1]{
\orgdiv{Centro de Matem\'atica e Aplica\c{c}\~oes, Departamento de Matem\'atica}, 
\orgname{Faculdade de Ci\^encias e Tecnologia, Universidade Nova de Lisboa}, 
\orgaddress{
\street{Quita da Torre}, 
\city{Caparica}, 
\postcode{2829--516}, 
\country{Portugal}}}
\affil[2]{
\orgdiv{Department of Mathematics}, 
\orgname{King's College London}, 
\orgaddress{
\street{Strand}, 
\city{London}, 
\postcode{WC2R 2LS},  
\country{United Kingdom}}}
\abstract{Let $X(\mathbb{R})$ be a separable translation-invariant Banach
function space and $a$ be a Fourier multiplier on $X(\mathbb{R})$. We prove
that the Wiener-Hopf operator $W(a)$ with symbol $a$ is maximally noncompact
on the space $X(\mathbb{R}_+)$, that is, its Hausdorff measure of 
noncompactness, its essential norm and its norm are all equal. This equality
for the Hausdorff measure of noncompactness of $W(a)$ 
is new even in the case of $X(\mathbb{R})=L^p(\mathbb{R})$ with $1\le p<\infty$.
}
\maketitle
\section{Introduction}
For Banach spaces $\mathcal{X}, \mathcal{Y}$, let 
$\mathcal{B}(\mathcal{X},\mathcal{Y})$ and 
$\mathcal{K}(\mathcal{X},\mathcal{Y})$ denote the sets of bounded linear and 
compact linear operators from $\mathcal{X}$ to $\mathcal{Y}$, 
respectively. We will abbreviate 
$\mathcal{B}(\mathcal{X}):=\mathcal{B}(\mathcal{X},\mathcal{X})$
and $\mathcal{K}(\mathcal{X}):=\mathcal{K}(\mathcal{X},\mathcal{X})$.
The norm of an operator $A\in\mathcal{B}(\mathcal{X},\mathcal{Y})$ is denoted 
by $\|A\|_{\mathcal{B}(\mathcal{X},\mathcal{Y})}$. The essential norm of 
$A \in \mathcal{B}(\mathcal{X},\mathcal{Y})$ is defined by
\[
\|A\|_{\mathcal{B}(\mathcal{X},\mathcal{Y}),\mathrm{e}} 
:= 
\inf\{\|A - K\|_{\mathcal{B}(\mathcal{X},\mathcal{Y})}\ : \  
K \in \mathcal{K}(\mathcal{X},\mathcal{Y})\} .
\]
For a bounded subset $\Omega$ of the space $\mathcal{X}$, we denote by 
$\chi(\Omega)$ the greatest lower bound of the set of numbers $r$ such that 
$\Omega$ can be covered by a finite family of open balls of radius $r$.
For $A \in \mathcal{B}(\mathcal{X},\mathcal{Y})$, set
\[
\|A\|_{\mathcal{B}(\mathcal{X},\mathcal{Y}),\chi} := 
\chi\left(A(B_\mathcal{X})\right) ,
\]  
where $B_\mathcal{X}$ denotes the closed unit ball in $\mathcal{X}$. The 
quantity $\|A\|_{\mathcal{B}(\mathcal{X},\mathcal{Y}),\chi}$ is 
called the Hausdorff measure of noncompactness of the operator $A$. It 
follows from the definition of the essential norm and 
\cite[inequality (3.29)]{LS71} that for every 
$A\in\mathcal{B}(\mathcal{X},\mathcal{Y})$ one has
\begin{equation}\label{eq:measure-essential-norm-arbitrary-operator}
\|A\|_{\mathcal{B}(\mathcal{X},\mathcal{Y}),\chi} 
\le 
\|A\|_{\mathcal{B}(\mathcal{X},\mathcal{Y}),\mathrm{e}}
\le
\|A\|_{\mathcal{B}(\mathcal{X},\mathcal{Y})}.
\end{equation}
Note that both inequalities in \eqref{eq:measure-essential-norm-arbitrary-operator}
may be strict. A simple example of an operator, for which the first inequality
is strict, is given in \cite[Section~2.4.11]{AKPRS92}. Moreover, there exist 
Banach spaces $\mathcal{X}$ and $\mathcal{Y}$ such that
\[
\sup_{A \in \mathcal{B}(\mathcal{X},\mathcal{Y})} 
\frac{\|A\|_{\mathcal{B}(\mathcal{X},\mathcal{Y}),\mathrm{e}}}
{\|A\|_{\mathcal{B}(\mathcal{X},\mathcal{Y}),\chi}} = \infty
\]
(see \cite[the proof of Theorem 2.5]{AT87}).

An interesting example of an operator, for which the second inequality in 
\eqref{eq:measure-essential-norm-arbitrary-operator} is strict, is the Cauchy 
singular integral operator
\[
(S_\Gamma f)(t):=\frac{1}{\pi i}\int_\Gamma \frac{f(\tau)}{\tau-t}d\tau
\]
on $L^2(\Gamma)$, where $\Gamma$ is a smooth simple closed contour 
in the complex plane. In this case, 
$\|S_\Gamma\|_{\mathcal{B}(L^2(\Gamma)),\mathrm{e}}=
\|S_\Gamma\|_{\mathcal{B}(L^2(\Gamma))}$
if and only if $\Gamma$ is a circle (see \cite[Corollary 5.7]{K10}).

We refer to the monographs \cite{AKPRS92,ADl97,BM14}
for the general theory of measures of noncompactness.
In particular, it is well known that
\[
K\in\mathcal{K}(\mathcal{X},\mathcal{Y})
\Longleftrightarrow
\|K\|_{\mathcal{B}(\mathcal{X},\mathcal{Y}),\chi}=0
\Longleftrightarrow
\|K\|_{\mathcal{B}(\mathcal{X},\mathcal{Y}),\mathrm{e}}=0
\]
(for the first equivalence, see, e.g., \cite[Theorem~5.29]{BM14}).
We will say that an operator $A\in\mathcal{B}(\mathcal{X},\mathcal{Y})$ 
is maximally noncompact if 
\[
\|A\|_{\mathcal{B}(\mathcal{X},\mathcal{Y}),\chi} =
\|A\|_{\mathcal{B}(\mathcal{X},\mathcal{Y}),\mathrm{e}}=
\|A\|_{\mathcal{B}(\mathcal{X},\mathcal{Y})}.
\]

Let $a\in L^\infty(\mathbb{R})$ and $X(\mathbb{R})$ be a Banach function space
(see \cite{L55}, Subsection~\ref{subsec:BFS} below, and also 
\cite[Ch.~1]{BS88}, \cite[Ch.~6]{PKJF13}). For a set $Q\subset\mathbb{R}$
of positive measure, let $e_Q$ be the the operator of 
extension by $0$ from $Q$ to $\mathbb{R}$ and $X(Q)$
be the Banach space of measurable functions $f:Q\to\mathbb{C}$ such that 
$\|f\|_{X(Q)}:=\|e_Qf\|_{X(\mathbb{R})}<\infty$.

Before discussing Wiener-Hopf operators, let us consider the multiplication operator 
$aI \in\mathcal{B}(X(\mathbb{R}))$. This operator is maximally noncompact. 
Indeed, it is easy to see that 
$\|aI\|_{\mathcal{B}(X(\mathbb{R}))} \le \|a\|_{L^\infty}$. On the 
other hand, for any $\varepsilon > 0$, there exist $c \in \mathbb{C}$ and a set 
of positive measure $Q \subset \mathbb{R}$ such that 
$|c| > \|a\|_{L^\infty} - \varepsilon/2$
and $|a - c| < \varepsilon/2$ a.e. in $Q$. Then
\begin{align*}
\|aI\|_{\mathcal{B}(X(\mathbb{R})),\chi} 
&\ge 
\|aI\|_{\mathcal{B}(X(Q)),\chi} 
\ge 
\|cI\|_{\mathcal{B}(X(Q)),\chi} - \|(c - a)I\|_{\mathcal{B}(X(Q)),\chi} 
\\
&\ge 
|c| \|I\|_{\mathcal{B}(X(Q)),\chi} - \|(c - a)I\|_{\mathcal{B}(X(Q))} 
> 
|c| - \varepsilon/2 > \|a\|_{L^\infty} - \varepsilon
\end{align*}
(see \cite[Theorem~1.1.6]{AKPRS92}). Since $\varepsilon > 0$ is arbitrary, one gets 
$\|aI\|_{\mathcal{B}(X(\mathbb{R})),\chi} \ge \|a\|_{L^\infty}$.
So, 
\[
\|aI\|_{\mathcal{B}(X(\mathbb{R})),\chi} = 
\|a\|_{L^\infty} = 
\|aI\|_{\mathcal{B}(X(\mathbb{R}))}.
\]

It is instructive to compare the above with what is known for multiplication 
operators on Sobolev spaces $W^m_p(\mathbb{R})$, $m \in \mathbb{N}$, 
$1 < p < \infty$. Suppose, for simplicity, that the support of $a$ is a subset 
of $[-1, 1]$. Then $\|aI\|_{\mathcal{B}(W^m_p(\mathbb{R}))}$ is estimated 
below by $\|a\|_{W^m_p}$ times a constant 
independent of $a$ (see \cite[Corollary~2.3.1]{MS09}).
On the other hand, $\|aI\|_{\mathcal{B}(W^m_p(\mathbb{R})),\chi}$ is 
estimated above by $\|a\|_{L^\infty}$ times a constant independent of 
$a$ (see \cite{ES05}). Hence, one cannot expect 
$aI : W^m_p(\mathbb{R}) \to W^m_p(\mathbb{R})$ to be maximally 
noncompact.

Let $S(\mathbb{R})$ denote the Schwartz spaces of rapidly
decaying infinitely differentiable functions, let
\[
(Fu)(\xi):=\widehat{u}(\xi):=\int_{\mathbb{R}}{u(x)}e^{-ix\xi}\,dx,
\quad
\xi\in\mathbb{R},
\]
be the Fourier transform of $u\in S(\mathbb{R})$, and let $F^{-1}$
denote the inverse Fourier transform on $S(\mathbb{R})$. The operators 
$F$ and $F^{-1}$ can be extended by continuity to $L^2(\mathbb{R})$, 
and we will use the same symbols for their extensions to 
$L^2(\mathbb{R})$. Let $X(\mathbb{R})$ be a Banach function space.
A function $a\in L^\infty(\mathbb{R})$ is said to belong to the set
$\mathcal{M}^0_{X(\mathbb{R})}$ of Fourier multipliers on 
$X(\mathbb{R})$ if
\[
\|a\|_{\mathcal{M}^0_{X(\mathbb{R})}}
:=
\sup\left\{
\frac{\|F^{-1} aFu\|_{X(\mathbb{R})}}{\|u\|_{X(\mathbb{R})}}:
u \in \left(L^2(\mathbb{R})\cap X(\mathbb{R})\right) \setminus\{0\}
\right\}<\infty.
\]
The operator $F^{-1} aF$ is translation-invariant, and one can find results on maximal noncompactness of such operators in \cite[Section 5]{KS24}.
If $a \in \mathcal{M}^0_{X(\mathbb{R})}$ and $X(\mathbb{R})$ is separable, in which case $L^2(\mathbb{R})\cap X(\mathbb{R})$ is dense 
in $X(\mathbb{R})$ (see Lemma \ref{le:density-simple} below), then $F^{-1} aF$ can be extended by continuity to a bounded linear operator
on $X(\mathbb{R})$. This operator is, in a sense, equivalent to the multiplication
operator $aI :  F(X(\mathbb{R})) \to F(X(\mathbb{R}))$, but maximal noncompactness of the former 
does not seem to follow directly from the maximal noncompactness of the operator $aI$ on the Banach function space $X(\mathbb{R})$.

Denote by $r_+$ the operator of restriction from $\mathbb{R}$ 
to $\mathbb{R}_+:=(0,\infty)$ and by $e_+$ the operator of 
extension by $0$ from $\mathbb{R}_+$ to $\mathbb{R}$. For a 
closed subspace $Y(\mathbb{R})$ of a Banach function space 
$X(\mathbb{R})$, let $Y(\mathbb{R}_+)$ be the space of all measurable 
functions $f:\mathbb{R}_+\to\mathbb{C}$ such that
\[
\|f\|_{Y(\mathbb{R}_+)}:=\|e_+ f\|_{Y(\mathbb{R})}<\infty.
\]
For $a \in \mathcal{M}^0_{X(\mathbb{R})}$, let
\[
W(a) u := r_+ F^{-1} aFe_+ u, 
\quad u \in L^2(\mathbb{R}_+)\cap X(\mathbb{R}_+) 
\]
be the Wiener-Hopf operator with symbol $a$. Put
\[
\|W(a)\|_{[X(\mathbb{R}_+)]} := 
\sup\left\{
\frac{\|W(a) u\|_{X(\mathbb{R}_+)}}{\|u\|_{X(\mathbb{R}_+)}}:
u \in \left(L^2(\mathbb{R}_+)\cap X(\mathbb{R}_+)\right) \setminus\{0\}
\right\} .
\]
Let $X_2(\mathbb{R}_+)$ be the closure of 
$L^2(\mathbb{R}_+)\cap X(\mathbb{R}_+)$ in $X(\mathbb{R}_+)$.
Then the bounded linear operator 
$W(a):L^2(\mathbb{R}_+)\cap X(\mathbb{R}_+)\to X(\mathbb{R}_+)$
can be extended by continuity to a bounded linear operator from 
$X_2(\mathbb{R}_+)$ to $X(\mathbb{R}_+)$, which we will denote again 
by $W(a)$. Moreover,
\[
\|W(a)\|_{[X(\mathbb{R}_+)]}
=
\|W(a)\|_{\mathcal{B}(X_2(\mathbb{R}_+),X(\mathbb{R}_+))}.
\]

We are interested in the maximal noncompactness of Wiener-Hopf operators.
We refer to \cite{KS22,KS24,EL25-1,EL25-2} for maximal noncompactness of some other 
integral transforms. 
 
Let us recall the results motivating our work.
If $a\in\mathcal{M}_{L^p(\mathbb{R})}^0$, then
\begin{equation}\label{eq:WH-norn-essential-norm}
\|W(a)\|_{\mathcal{B}(L^p(\mathbb{R}_+))}
=
\|W(a)\|_{\mathcal{B}(L^p(\mathbb{R}_+)),\mathrm{e}},
\quad
1<p<\infty.
\end{equation}
(see \cite[p.~5]{BGS97}, where it is stated without proof
for $1\le p<\infty$, and \cite[Section~9.5(b)]{BS06},
where the exact range for $p$ is not specified and the 
given idea of the proof seems to exclude the case $p=1$).

Equality \eqref{eq:WH-norn-essential-norm}
was recently extended by the first author and M.~Valente
to the setting of Lorentz and Orlicz spaces
(see, e.g., \cite[Ch.~4]{BS88} for their definitions).
Let $X(\mathbb{R})$ be a Lorentz space $L^{p,q}(\mathbb{R})$
with $1<p<\infty$ and $1\le q<\infty$ or a separable Orlicz space
$L^\Phi(\mathbb{R})$ generated by a Young's function $\Phi$ satisfying
\begin{equation}\label{eq:Orlicz-excluding-L1}
\lim_{t\to 0^+}\frac{\Phi(t)}{t}=0.
\end{equation}
It was shown in \cite[Theorem~1]{KV25} that if 
$a\in\mathcal{M}_{X(\mathbb{R})}^0$, then
\begin{equation}
\|W(a)\|_{\mathcal{B}(X(\mathbb{R_+}))}
=
\|W(a)\|_{\mathcal{B}(X(\mathbb{R_+})),\mathrm{e}}.
\label{eq:Karlovych-Valente}
\end{equation}
It is well-known that the Lebesque space $L^p(\mathbb{R})$, $1\le p<\infty$, 
is the Orlicz space generated by the Young's function $\Phi(t)=t^p$
(with equal norms). So, condition \eqref{eq:Orlicz-excluding-L1} excludes 
the case of the space $L^1(\mathbb{R})$. 

Note that Lebesgue, Lorentz, and Orlicz spaces are translation-invariant.
So it is natural to look for extensions of the above results to the setting
of arbitrary translation-invariant Banach function spaces (see 
Subsection~\ref{subsec:TI-BFS} for their definition). Let us state here a
simply formulated consequence of our main result, which will be proved 
in Subsection~\ref{subsec:main-result}.
\begin{theorem}\label{th:nice}
Suppose $X(\mathbb{R})$ is a separable translation-invariant Banach function
space. If $a\in\mathcal{M}_{X(\mathbb{R})}^0$, then the Wiener-Hopf operator 
$W(a)$ is maximally noncompact on the space $X(\mathbb{R}_+)$, that is,
\[
\|W(a)\|_{\mathcal{B}(X(\mathbb{R}_+))}
=
\|W(a)\|_{\mathcal{B}(X(\mathbb{R}_+)),\mathrm{e}}
=
\|W(a)\|_{\mathcal{B}(X(\mathbb{R}_+)),\chi}.
\]
\end{theorem}
This theorem improves equality \eqref{eq:Karlovych-Valente} for separable
Orlicz spaces because it allows to drop condition 
\eqref{eq:Orlicz-excluding-L1} (and, so, to include the space $L^1(\mathbb{R}_+)$
into consideration). What is more important, is that the equality 
$\|W(a)\|_{\mathcal{B}(X(\mathbb{R}_+))}
=
\|W(a)\|_{\mathcal{B}(X(\mathbb{R}_+)),\chi}$
is new even for Lebesgue spaces $L^p(\mathbb{R}_+)$ with $1\le p<\infty$.

The paper is organised as follows. In Section~\ref{sec:preliminaries}, we collect
preliminaries on the class of Banach function spaces and its subclasses of
translation-invariant and rear\-range\-ment-in\-vari\-ant spaces. For the latter
subclass of spaces, we recall the definition of their Zippin (or fundamental)
indices. Section~\ref{sec:main-result-and-consequences} starts with our
main result (Theorem~\ref{th:main-result}) on the maximal noncompactness
of the truncation $A_+:=r_+Ae_+$ of a translation-invariant operator
$A:X(\mathbb{R})\to Y(\mathbb{R})$  acting from a translation-invariant
subspace $X(\mathbb{R})$ to a translation-invariant Banach function space
$Y(\mathbb{R})$ with some additional properties. Further, this result is applied
to prove Theorem~\ref{th:nice}. Finally, we show that the hypotheses
of Theorem~\ref{th:main-result} are flexible enough to prove analogues
of Theorem~\ref{th:nice} for a not necessarily separable rearrangement-invariant
Banach functionn space with the positive lower Zippin index. In particular, we state
such a result for weak Lebesgue spaces $L^{p,\infty}$ with $1<p<\infty$, which 
are not separable, as it is well known.
\section{Preliminaries on Banach function spaces}\label{sec:preliminaries}
\subsection{Banach function spaces}\label{subsec:BFS}
The set of all Lebesgue measurable complex-valued functions on $\mathbb{R}$
is denoted by $\mathfrak{M}(\mathbb{R})$. Let $\mathfrak{M}^+(\mathbb{R})$
be the subset of functions  in $\mathfrak{M}(\mathbb{R})$ whose values lie
in $[0,\infty]$. The characteristic function of a measurable set
$E\subset\mathbb{R}$ is denoted by $\mathbbm{1}_E$ and the Lebesgue measure 
of $E$ is denoted by $|E|$.

Following \cite[p.~3]{L55} (see also \cite[Ch.~1, Definition~1.1]{BS88} and 
\cite[Definition~6.1.5]{PKJF13}), a mapping
$\rho:\mathfrak{M}^+(\mathbb{R})\to [0,\infty]$
is called a Banach function norm if, for all functions $f,g,f_j$
($j\in\mathbb{N}$) in $\mathfrak{M}^+(\mathbb{R})$, for all constants
$a\ge 0$, and for all measurable subsets $E$ of $\mathbb{R}$, the
following properties hold:
\begin{eqnarray*}
{\rm (A1)} &\quad &
\rho(f)=0  \Leftrightarrow  f=0\ \mbox{a.e.},
\quad
\rho(af)=a\rho(f),
\quad
\rho(f+g) \le \rho(f)+\rho(g),\\
{\rm (A2)} &\quad &
0\le g \le f \ \mbox{a.e.} \ \Rightarrow \ \rho(g)
\le \rho(f)
\quad\mbox{(the lattice property)},\\
{\rm (A3)} &\quad &
0\le f_j \uparrow f \ \mbox{a.e.} \ \Rightarrow \
       \rho(f_j) \uparrow \rho(f)\quad\mbox{(the Fatou property)},\\
{\rm (A4)} &\quad &
E \text{ is bounded} \Rightarrow \rho(\mathbbm{1}_E) <\infty,\\
{\rm (A5)} &\quad &
E \text{ is bounded} \Rightarrow \int_E f(x)\,dx \le C_E\rho(f)
\end{eqnarray*}
with $C_E \in (0,\infty)$ that may depend on $E$ and $\rho$ but is
independent of $f$.

When functions differing only on a set of measure zero are identified, the
set $X(\mathbb{R})$ of all functions $f\in\mathfrak{M}(\mathbb{R})$ for
which $\rho(|f|)<\infty$ becomes a Banach space under the norm
\[
\|f\|_{X(\mathbb{R})} :=\rho(|f|)
\]
and  under the natural linear space operations
(see \cite[Ch.~1, \S~1, Theorem~1]{L55} or
\cite[Ch.~1, Theorems~1.4 and~1.6]{BS88}). It is  called a 
{\it Banach function space}. 

If $\rho$ is a Banach function norm, its associate norm $\rho'$ is defined on
$\mathfrak{M}^+(\mathbb{R})$ by
\[
\rho'(g):=\sup\left\{
\int_{\mathbb{R}} f(x)g(x)\,dx \ : \
f\in \mathfrak{M}^+(\mathbb{R}), \ \rho(f) \le 1
\right\}.
\]
It is a Banach function norm itself (see \cite[Ch.~1, \S~1]{L55} or
\cite[Ch.~1, Theorem~2.2]{BS88}).
The Banach function space $X'(\mathbb{R})$ defined by the Banach
function norm $\rho'$ is called the associate space (K\"othe dual) of
$X(\mathbb{R})$. The Lebesgue space $L^p(\mathbb{R})$, $1\le p\le\infty$,
is the archetypical example of a Banach function space. Other classical
examples of Banach function spaces are 
Lorentz spaces $L^{p,q}(\mathbb{R})$ with $1<p<\infty$ and $1\le q\le\infty$ 
(see, e.g., \cite[Ch.~4, Section~4]{BS88}),
Orlicz spaces $L^\Phi(\mathbb{R})$
(see, e.g., \cite[Ch.~4, Section~8]{BS88}), 
all other rearrangement-invariant Banach function spaces
 (see \cite[Ch.~2]{BS88}), 
as well as, variable Lebesgue spaces $L^{p(\cdot)}(\mathbb{R})$ 
(see \cite{CF13}).
\begin{remark}
We note that our definition of a Banach function space is slightly 
different from that found in \cite[Ch.~1, Definition~1.1]{BS88}
and \cite[Definition~6.1.5]{PKJF13}. 
In particular, in Axioms (A4) and (A5)
we assume that the set $E$ is a bounded set, whereas it is sometimes 
assumed that $E$ merely satisfies $|E| <\infty$.  It is well known that
all main elements of the general theory of Banach function spaces
work with (A4) and (A5) stated for bounded sets \cite{L55} (see also 
the discussion at the beginning of Chapter~1 on page~2 of \cite{BS88} 
and \cite{LN24}).
\end{remark}
\subsection{Density of simple functions and its consequences}
Let $S_0(\mathbb{R})$ denote the set of all simple functions with compact
support in $\mathbb{R}$ and $S_0(\mathbb{R}_+)$ denote the set of all
simple functions with compact support in $[0,\infty)$.
\begin{lemma}\label{le:density-simple}
If $X(\mathbb{R})$ is a separable Banach function space, then the sets
$S_0(\mathbb{R})$, $L^2(\mathbb{R})\cap X(\mathbb{R})$ are dense
in the space $X(\mathbb{R})$, and the sets $S_0(\mathbb{R}_+)$,
 $L^2(\mathbb{R}_+)\cap X(\mathbb{R}_+)$ are dense in the
 space $X(\mathbb{R}_+)$.
\end{lemma}
\begin{proof}
Since $S_0(\mathbb{R})\subset L^2(\mathbb{R})\cap X(\mathbb{R})$
and $S_0(\mathbb{R}_+)\subset L^2(\mathbb{R}_+)\cap X(\mathbb{R}_+)$,
it is enough to prove that $S_0(\mathbb{R})$ is dense in $X(\mathbb{R})$
and $S_0(\mathbb{R}_+)$ is dense in $X(\mathbb{R}_+)$.

Following \cite[Ch.~1, \S~2, Definition~3]{L55}, let $X_b(\mathbb{R})$
be the closure of the set of all bounded measurable functions with compact support.
It can be shown as in \cite[Ch.~1, Proposition~3.10]{BS88} that
the closure of $S_0(\mathbb{R})$ in $X(\mathbb{R})$ coincides with
$X_b(\mathbb{R})$. It follows from \cite[Ch.~1, \S~2, Lemma~4 and \S~3, Corollary~1]{L55}
(see alos \cite[Ch.~1, Theorem~3.11 and Corollary~5.6]{BS88}) that if
$X(\mathbb{R})$ is separable, then $X(\mathbb{R})=X_b(\mathbb{R})$.
So, in this case $S_0(\mathbb{R})$ is dense in $X(\mathbb{R})$.

Now, let $f\in X(\mathbb{R}_+)$. Then $g:=e_+f\in X(\mathbb{R})$
and for every $\varepsilon>0$ there exists $h\in S_0(\mathbb{R})$
such that $\|g-h\|_{X(\mathbb{R})}<\varepsilon$. Therefore
$r:=r_+h\in S_0(\mathbb{R}_+)$ and
\begin{align*}
\|f-r\|_{X(\mathbb{R}_+)}
&=
\|r_+e_+f-r_+h\|_{X(\mathbb{R}_+)}
=
\|r_+(g-h)\|_{X(\mathbb{R}_+)}
\\
&=
\|e_+r_+(g-h)\|_{X(\mathbb{R})}
=
\|\mathbbm{1}_{(0,\infty)}(g-h)\|_{X(\mathbb{R})}
\le 
\|g-h\|_{X(\mathbb{R})}<\varepsilon.
\end{align*}
This implies that $S_0(\mathbb{R}_+)$ is dense in $X(\mathbb{R}_+)$.
\end{proof}
Let us prove another result where an argument based on approximation
of functions in $X'(\mathbb{R})$ by functions in $S_0(\mathbb{R})$
plays an important role (see the proof of \cite[Lemma~2.10]{KS19}).
\begin{lemma}\label{le:norm-via-simple-functions}
Let $X(\mathbb{R})$ be a Banach function space and let $X'(\mathbb{R})$ 
be its associate space. For every $f\in X(\mathbb{R}_+)$, one has
\begin{equation}\label{eq:norm-via-simple-functions-1}
\|f\|_{X(\mathbb{R}_+)}
=
\sup\left\{
\left|\int_{\mathbb{R}_+}f(x)s_+(x)\,dx\right| \ :\ 
s_+\in S_0(\mathbb{R}_+), \
\|s_+\|_{X'(\mathbb{R}_+)}\le 1
\right\}.
\end{equation}
\end{lemma}
\begin{proof}
It follows from \cite[Lemma~2.10]{KS19} that
\begin{equation}\label{eq:norm-via-simple-functions-2}
\|f\|_{X(\mathbb{R}_+)}
=
\|e_+f\|_{X(\mathbb{R})}
=
\sup\left\{
\left|\int_{\mathbb{R}}(e_+f)(x)s(x)\,dx\right| \ :\ 
s\in S_0(\mathbb{R}), \
\|s\|_{X'(\mathbb{R})}\le 1
\right\}.
\end{equation}
If $s\in S_0(\mathbb{R})$ and $\|s\|_{X'(\mathbb{R})}\le 1$, then 
$s_+=r_+s \in S_0(\mathbb{R}_+)$,
\[
\|s_+\|_{X'(\mathbb{R}_+)}
=
\|e_+r_+s\|_{X'(\mathbb{R})}
=
\|\mathbbm{1}_{(0,\infty)}s\|_{X'(\mathbb{R})}
\le
\|s\|_{X'(\mathbb{R})}
\le 1,
\]
and
\begin{equation}\label{eq:norm-via-simple-functions-3}
\left|\int_{\mathbb{R}}(e_+f)(x)s(x)\,dx\right| 
=
\left|\int_{\mathbb{R}_+}f(x)s_+(x)\,dx\right|.
\end{equation}
Hence
\begin{align}
&
\sup\left\{
\left|\int_{\mathbb{R}}(e_+f)(x)s(x)\,dx\right| \ :\ 
s\in S_0(\mathbb{R}), \
\|s\|_{X'(\mathbb{R})}\le 1
\right\}
\nonumber\\
&\quad \le 
\sup\left\{
\left|\int_{\mathbb{R}_+}f(x)s_+(x)\,dx\right| \ :\ 
s_+\in S_0(\mathbb{R}_+), \
\|s_+\|_{X'(\mathbb{R}_+)}\le 1
\right\}.
\label{eq:norm-via-simple-functions-4}
\end{align}
On the other hand, if $s_+\in S_0(\mathbb{R}_+)$ and
$\|s_+\|_{X'(\mathbb{R}_+)}\le 1$, then $s=e_+s_+\in S_0(\mathbb{R})$ and
\[
\|s\|_{X'(\mathbb{R})}=\|e_+s_+\|_{X'(\mathbb{R})}=\|s_+\|_{X'(\mathbb{R}_+)}\le 1,
\]
and \eqref{eq:norm-via-simple-functions-3} holds. Therefore
\begin{align}
&
\sup\left\{
\left|\int_{\mathbb{R}_+}f(x)s_+(x)\,dx\right| \ :\ 
s_+\in S_0(\mathbb{R}_+), \
\|s_+\|_{X'(\mathbb{R}_+)}\le 1
\right\}
\nonumber\\
&\quad \le 
\sup\left\{
\left|\int_{\mathbb{R}}(e_+f)(x)s(x)\,dx\right| \ :\ 
s\in S_0(\mathbb{R}), \
\|s\|_{X'(\mathbb{R})}\le 1
\right\}.
\label{eq:norm-via-simple-functions-5}
\end{align}
Combining \eqref{eq:norm-via-simple-functions-2}
and \eqref{eq:norm-via-simple-functions-4}--\eqref{eq:norm-via-simple-functions-5},
we arrive at \eqref{eq:norm-via-simple-functions-1}.
\end{proof}
\subsection{Translation-invariant subspaces of Banach function spaces}
\label{subsec:TI-BFS}
A closed subspace $Y(\mathbb{R})$ of a Banach function space
$X(\mathbb{R})$ is said to by translation-invariant if for all 
$y \in \mathbb{R}$ and for all functions $u \in Y(\mathbb{R})$, one has 
$\tau_y u\in Y(\mathbb{R})$ and
\[
\|\tau_y u\|_{Y(\mathbb{R})} =
\|\tau_y u\|_{X(\mathbb{R})} = 
\|u\|_{X(\mathbb{R})}= 
\|u\|_{Y(\mathbb{R})}, 
\]
where the translation operator $\tau_y$ is defined by
$(\tau_y u)(x) := u(x - y)$ for all $x \in \mathbb{R}$.
Here and in what follows we always suppose that
subspaces of Banach function spaces 
are equipped with the induced norms.

Note that all rearrangement-invariant Banach function spaces 
(see \cite[Ch.~2]{BS88} or Subsction~\ref{subsec:RI-BFS} below) are
translation-invariant. On the other hand, in view of 
\cite[Theorem~5.17]{CF13}, variable Lebesgue spaces 
$L^{p(\cdot)}(\mathbb{R})$ are not translation-invariant.

We will repeatedly use the fact that a Banach function space $X(\mathbb{R})$
is translation-invariant if and only if its associate space $X'(\mathbb{R})$
is translation-invariant (see \cite[Lemma~2.1]{KS19}). 
\subsection{Rearrangement-invariant Banach function spaces and their Zippin indices}
\label{subsec:RI-BFS}
Let $\mathfrak{M}_0(\mathbb{R})$ and  $\mathfrak{M}_0^+(\mathbb{R})$ 
be the classes of a.e. finite functions in $\mathfrak{M}(\mathbb{R})$ and 
$\mathfrak{M}^+(\mathbb{R})$, respectively. The distribution function 
$\mu_f$ of $f\in\mathfrak{M}_0(\mathbb{R})$ is given by
\[
\mu_f(\lambda):=
|\{x\in\mathbb{R}:|f(x)|>\lambda\}|,
\quad
\lambda\ge 0.
\]
The non-increasing rearrangement of $f\in\mathfrak{M}_0(\mathbb{R})$ is 
the function defined by
\[
f^*(t):=\inf\{\lambda:\mu_f(\lambda)\le t\},\quad t\ge 0.
\]
We use here the standard convention that $\inf\emptyset=+\infty$. 
Two functions $f\in\mathfrak{M}_0(\mathbb{R})$ 
and $g\in\mathfrak{M}_0(\mathbb{R})$
are said to be equimeasurable if 
$\mu_f(\lambda)=\mu_g(\lambda)$ for all $\lambda\ge 0$.

A Banach function norm $\rho:\mathfrak{M}^+(\mathbb{R}) \to [0,\infty]$ 
is called rearrangement-in\-vari\-ant if for every pair of 
equimeasurable functions $f,g \in \mathfrak{M}^+_0(\mathbb{R})$, the 
equality $\rho(f)=\rho(g)$ holds. In that case, the Banach function space 
$X(\mathbb{R})$ generated by $\rho$ is said to be a 
rearrangement-invariant Banach function space (or simply a
rearrangement-invariant space). Le\-besgue spaces $L^p(\mathbb{R})$,
$1\le p\le\infty$, Orlicz spaces $L^\Phi(\mathbb{R})$, and Lorentz spaces 
$L^{p,q}(\mathbb{R})$ are  classical examples of rearrangement-inva\-riant 
Banach function spaces  (see, e.g., \cite{BS88} and the references therein). 
By  \cite[Ch.~2, Proposition~4.2]{BS88}, if a Banach function space 
$X(\mathbb{R})$ is rearrangement-invariant, then its associate space 
$X'(\mathbb{R})$ is also rearrangement-invariant.

Let $X(\mathbb{R})$ be a rearrangement-invariant Banach function space.
For each $t\in[0,\infty)$, let $E$ be a subset of $\mathbb{R}$ with
$|E|=t$ and let $\varphi_X(t)=\|\mathbbm{1}_E\|_{X(\mathbb{R})}$.
The function $\varphi_X$ so defined is called the fundamental function
of $X(\mathbb{R})$. Consider the function
\[
M_X(t):=\sup_{0<s<\infty}\frac{\varphi_X(st)}{\varphi_X(s)},
\quad
0<t<\infty.
\]
The numbers
\[
p_X
:=
\sup_{0<t<1}\frac{\ln M_X(t)}{\ln t}
=
\lim_{t\to 0^+}\frac{\ln M_X(t)}{\ln t},
\quad
q_X
:=
\inf_{1<t<\infty}\frac{\ln M_X(t)}{\ln t}
=
\lim_{t\to \infty}\frac{\ln M_X(t)}{\ln t}
\]
are called the lower and upper Zippin (or fundamental) indices of the space $X(\mathbb{R})$.
It is well known that $0\le p_X\le q_X\le 1$  
(see \cite[Exercise~14 to Chapter~3]{BS88} and \cite[Section~4]{M85}).  
Simple calculations show that $\varphi_{L^{p}}(t)=t ^{1/p}$ and hence 
$p_{L^p}=q_{L^p}=1/p$ for all $1\le p\le\infty$.
\section{Main result and its applications to Wiener-Hopf operators}\label{sec:main-result-and-consequences}
\subsection{Main result}\label{subsec:main-result}
Let $X(\mathbb{R})$ and $Y(\mathbb{R})$ be translation-invariant subspaces
of Banach function spaces. An operator 
$A\in\mathcal{B}(X(\mathbb{R}),Y(\mathbb{R}))$ is said to be translation
invariant if $\tau_y A=A\tau_y$ for every $y\in\mathbb{R}$.
For an operator $A : X(\mathbb{R}) \to Y(\mathbb{R})$, let
\[
A_+ u := r_+  Ae_+ u , \quad u \in  X(\mathbb{R}_+) . 
\]
For $n\in\mathbb{N}$, define
\[
V_n:=r_+\tau_n e_+.
\]
If $X(\mathbb{R})$ is a translation-invariant subspace of a Banach function space, then
for every $f\in X(\mathbb{R}_+)$ and $n\in\mathbb{N}$, one has
$\operatorname{supp}\tau_ne_+f\subset[n,\infty)$ and
\begin{align}
\|V_nf\|_{X(\mathbb{R}_+)}
&=
\|e_+r_+\tau_n e_+f\|_{X(\mathbb{R})}
=
\|\mathbbm{1}_{(0,\infty)}\tau_ne_+f\|_{X(\mathbb{R})}
\nonumber\\
&=
\|\tau_ne_+f\|_{X(\mathbb{R})}
=
\|e_+f\|_{X(\mathbb{R})}
=
\|f\|_{X(\mathbb{R}_+)}.
\label{eq:Vn-contraction}
\end{align}

If a space $Y(\mathbb{R}_+)$ is non-separable, then 
$S_0(\mathbb{R}_+)$ is not dense in it, but it may happen that every element 
of $Y(\mathbb{R}_+)$ can be approximated by elements of $S_0(\mathbb{R}_+)$ 
in a norm weaker than $\|\cdot\|_{Y(\mathbb{R}_+)}$. This possibility is described 
in the next theorem by introducing an auxiliary 
Banach function space $Z(\mathbb{R})$.
Note that the norm $\|\cdot\|_{Y(\mathbb{R}_+)}$
is, in general, different from the restriction of the norm 
$\|\cdot\|_{Z(\mathbb{R}_+)}$ to $Y(\mathbb{R}_+)$.
\begin{theorem}\label{th:main-result} 
Let $X(\mathbb{R})$ be a translation-invariant subspace of a Banach function 
space,  $Y(\mathbb{R})$ and $Z(\mathbb{R})$  be translation-invariant 
Banach function spaces such that $Y(\mathbb{R}_+)$ is a subset of 
the closure of $S_0(\mathbb{R}_+)$ in $Z(\mathbb{R}_+)$. If 
$A \in \mathcal{B}(X(\mathbb{R}), Y(\mathbb{R}))$
is a translation-invariant operator, then 
$A_+\in\mathcal{B}(X(\mathbb{R}_+),Y(\mathbb{R}_+))$ 
is maximally noncompact, that is,
\[
\|A_+\|_{\mathcal{B}(X(\mathbb{R}_+),Y(\mathbb{R}_+)),\chi} 
=  
\|A_+\|_{\mathcal{B}(X(\mathbb{R}_+),Y(\mathbb{R}_+)),\mathrm{e}} 
= 
\|A_+\|_{\mathcal{B}(X(\mathbb{R}_+),Y(\mathbb{R}_+))} .
\]
\end{theorem}

\begin{proof}
The proof is similar to the proofs of \cite[Theorems~1.1--1.2]{KS22}
and \cite[Theorem~1.1]{KS24}.
In view of \eqref{eq:measure-essential-norm-arbitrary-operator}, it is sufficient to prove that 
\begin{equation}\label{eq:main-proof-1}
\|A_+\|_{\mathcal{B}(X(\mathbb{R}_+),Y(\mathbb{R}_+)),\chi} 
\ge
\|A_+\|_{\mathcal{B}(X(\mathbb{R}_+),Y(\mathbb{R}_+))} .
\end{equation}
Take an arbitrary $\varepsilon > 0$. There exists $g \in X(\mathbb{R}_+)$ 
such that $\|g\|_{X(\mathbb{R}_+)} = 1$ 
and  
\begin{equation}\label{eq:main-proof-2}
\|A_+g\|_{Y(\mathbb{R}_+)} > \|A_+\|_{\mathcal{B}(X(\mathbb{R}_+),Y(\mathbb{R}_+))} - \varepsilon .
\end{equation}
It follows from Lemma~\ref{le:norm-via-simple-functions} that there exists 
$s \in S_0(\mathbb{R}_+)\setminus\{0\}$ such that $\|s\|_{Y'(\mathbb{R}_+)} \le 1$ and
\begin{equation}\label{eq:main-proof-3}
\left|\int_{\mathbb{R}_+} (A_+g)(x) s(x)\, dx \right| 
\ge 
\|A_+g\|_{Y(\mathbb{R}_+)} - \varepsilon.
\end{equation}
Set $s_n := V_n s$. Since $Y(\mathbb{R})$ is translation-invariant, 
in view of \cite[Lemma 2.1]{KS19}, so is $Y'(\mathbb{R})$. Therefore,
in view of \eqref{eq:Vn-contraction}, one has
\begin{equation}\label{eq:main-proof-4}
\|s_n\|_{Y'(\mathbb{R}_+)} = \|s\|_{Y'(\mathbb{R}_+)} \le 1,
\quad
n\in\mathbb{N}.
\end{equation}
Making a change of variables in the left-hand side of 
\eqref{eq:main-proof-3}, we see that for all $n\in\mathbb{N}$,
\[
\left|\int_{\mathbb{R}_+} (V_nA_+g)(x) s_n(x)\, dx \right| =
\left|\int_{\mathbb{R}_+} (A_+g)(x) s(x)\, dx \right| 
\ge 
\|A_+g\|_{Y(\mathbb{R}_+)} - \varepsilon .
\]
Since $A$ is translation-invariant, it is easy to see that
for all $n\in\mathbb{N}$,
\begin{align*}
V_nA_+g
&=
(r_+ \tau_n e_+)(r_+Ae_+) g
=
r_+ \tau_n \mathbbm{1}_{(0,\infty)} A e_+ g
=
r_+ \mathbbm{1}_{(n,\infty)} \tau_n A e_+ g
\\
&=
\mathbbm{1}_{(n,\infty)} r_+ A\tau_n e_+ g
=
\mathbbm{1}_{(n,\infty)} r_+ A \mathbbm{1}_{(0,\infty)} \tau_n e_+ g
=
\mathbbm{1}_{(n,\infty)} (r_+ A e_+)(r_+ \tau_n e_+) g
\\
&=
\mathbbm{1}_{(n,\infty)}A_+ V_n g.
\end{align*}
On the other hand, $\operatorname{supp} s_n \subset [n, \infty)$. Hence, 
$(A_+V_ng)s_n =  (V_nA_+g)s_n$ a.e., and one gets
\begin{equation}\label{eq:main-proof-5}
\left|\int_{\mathbb{R}_+} (A_+V_ng)(x) s_n(x)\, dx \right| \ge 
\|A_+g\|_{Y(\mathbb{R}_+)} - \varepsilon .
\end{equation}
Take any finite set $\{\varphi_1, \dots, \varphi_m\} \subset Y(\mathbb{R}_+)$. Since $Y(\mathbb{R}_+)$ is 
a subset of the closure of $S_0(\mathbb{R}_+)$ in $Z(\mathbb{R}_+)$, and 
$s\in S_0(\mathbb{R}_+)\subset Z'(\mathbb{R}_+)$, there exists a set
$\{\psi_1, \dots, \psi_m\} \subset S_0(\mathbb{R}_+)$ such that
\begin{equation}\label{eq:main-proof-6}
\|\varphi_j - \psi_j\|_{Z(\mathbb{R}_+)}  < \frac{\varepsilon}{\|s\|_{Z'(\mathbb{R}_+)}}\, , 
\quad 
j\in\{ 1, \dots, m\}.
\end{equation}
Taking into account that $\psi_j$ and $s$ are compactly supported, 
we see that there exists 
$N\in\mathbb{N}$ such that
\begin{equation}\label{eq:main-proof-7}
\psi_j s_N = 0,
\quad 
j\in\{ 1, \dots, m\}.
\end{equation}
Since $Z'(\mathbb{R})$ is translation-invariant (see \cite[Lemma 2.1]{KS19}),
it follows from \eqref{eq:Vn-contraction} that
\[
\|s_N\|_{Z'(\mathbb{R}_+)} = \|s\|_{Z'(\mathbb{R}_+)}. 
\]
Then, in view of 
equalities \eqref{eq:main-proof-7}, H\"older's inequality for the 
space $Z(\mathbb{R})$ (see \cite[Ch.~1, Theorem~2.4]{BS88}), 
and inequalities \eqref{eq:main-proof-6},
one gets, for $j\in\{1,\dots,m\}$,
\begin{align}
\left|\int_{\mathbb{R}_+} \varphi_j(x) s_N(x)\, dx\right| 
 & = 
\left|\int_{\mathbb{R}_+} 
\left(\varphi_j(x) - \psi_j(x)\right) s_N(x)\, dx
\right| 
\le 
\|\varphi_j - \psi_j\|_{Z(\mathbb{R}_+)} \|s_N\|_{Z'(\mathbb{R}_+)} 
\nonumber\\
&< 
\frac{\varepsilon}{\|s\|_{Z'(\mathbb{R}_+)}}\, \|s\|_{Z'(\mathbb{R}_+)} 
= 
\varepsilon.
\label{eq:main-proof-8}
\end{align}
Combining \eqref{eq:main-proof-5} and \eqref{eq:main-proof-8}, we see
that for $j\in\{1,\dots,m\}$,
\begin{align}
&
\left|
\int_{\mathbb{R}_+}  \big((A_+V_N g)(x) - \varphi_j(x)\big) s_N(x)\, dx
\right| 
\nonumber
\\
&\qquad\ge 
\left|\int_{\mathbb{R}_+} (A_+V_Ng)(x) s_N(x)\, dx \right|
- 
\left|\int_{\mathbb{R}_+} \varphi_j(x) s_N(x)\, dx\right| 
> 
\|A_+g\|_{Y(\mathbb{R}_+)} - 2\varepsilon .
\label{eq:main-proof-9}
\end{align}
On the other hand, applying H\"older's inequality to the space $Y(\mathbb{R})$
(see \cite[Ch.~1, Theorem~2.4]{BS88}) and taking into account
inequality \eqref{eq:main-proof-4}, we get for 
$j\in\{1,\dots,m\}$,
\begin{align}
\left|
\int_{\mathbb{R}_+}  \big((A_+V_N g)(x) - \varphi_j(x)\big) s_N(x)\, dx
\right| 
&\le 
\|A_+V_N g - \varphi_j\|_{Y(\mathbb{R}_+)} \|s_N\|_{Y'(\mathbb{R}_+)} 
\nonumber
\\
& \le 
\|A_+V_N g - \varphi_j\|_{Y(\mathbb{R}_+)} .
\label{eq:main-proof-10}
\end{align}
So, we deduce from \eqref{eq:main-proof-9}, \eqref{eq:main-proof-10} 
and \eqref{eq:main-proof-2} that for $j\in\{1,\dots,m\}$,
\begin{align*}
\|A_+V_N g - \varphi_j\|_{Y(\mathbb{R}_+)}
>  
\|A_+g\|_{Y(\mathbb{R}_+)} - 2\varepsilon 
> 
\|A_+\|_{\mathcal{B}(X(\mathbb{R}_+),Y(\mathbb{R}_+))}  - 3\varepsilon.
\end{align*}
Since $X(\mathbb{R})$ is  a translation-invariant subspace of a Banach function space, 
it follows from \eqref{eq:Vn-contraction} that 
\[
\|V_N g\|_{X(\mathbb{R}_+)} = \|g\|_{X(\mathbb{R}_+)} = 1.
\]
So, for every finite set 
$\{\varphi_1, \dots, \varphi_m\} \subset Y(\mathbb{R}_+)$, there exist an element 
$A_+V_Ng$ of the image of the unit ball $A_+\big(B_{X(\mathbb{R}_+)}\big)$ that lies 
at a distance greater than $\|A_+\|_{\mathcal{B}(X(\mathbb{R}_+),Y(\mathbb{R}_+))} - 3\varepsilon$ 
from every element of $\{\varphi_1, \dots, \varphi_m\}$.
This means that $A_+\big(B_{X(\mathbb{R}_+)}\big)$ cannot 
be covered by a finite family of open balls of radius 
$\|A_+\|_{\mathcal{B}(X(\mathbb{R}_+),Y(\mathbb{R}_+))}  - 3\varepsilon$. Hence, for all $\varepsilon>0$,
\[
\|A_+\|_{\mathcal{B}(X(\mathbb{R}_+),Y(\mathbb{R}_+)),\chi} 
\ge 
\|A_+\|_{\mathcal{B}(X(\mathbb{R}_+),Y(\mathbb{R}_+))}  - 3\varepsilon. 
\]
Passing to the limit as $\varepsilon\to 0+$, we 
arrive at \eqref{eq:main-proof-1}, which completes the proof.
\end{proof}
\subsection{Proof of Theorem~\ref{th:nice}}
If $X(\mathbb{R})$ is a translation-invariant Banach function space, then its subspace
$X_2(\mathbb{R})$, being the closure of $L^2(\mathbb{R})\cap X(\mathbb{R})$
in $X(\mathbb{R})$, is also translation-invariant. It is easy to see that
if $a\in\mathcal{M}_{X(\mathbb{R})}^0$, then the Fourier multiplier
operator 
\begin{equation}\label{eq:Fourier-convolution-operator}
A=F^{-1}aF:X_2(\mathbb{R})\to X(\mathbb{R}) 
\end{equation}
is translation-invariant and
$A_+=W(a)\in\mathcal{B}(X_2(\mathbb{R}_+),X(\mathbb{R}_+))$.

It follows from Lemma~\ref{le:density-simple} that if $X(\mathbb{R})$ is separable,
then $X_2(\mathbb{R})=X(\mathbb{R})$ and $S_0(\mathbb{R}_+)$ is dense
in $X(\mathbb{R}_+)$. So, applying Theorem~\ref{th:main-result} to $A$, $A_+$
as above and $X(\mathbb{R})=X_2(\mathbb{R})=Y(\mathbb{R})=Z(\mathbb{R})$,
we arrive at Theorem~\ref{th:nice}.
\qed
\subsection{Maximal noncompactness of Wiener-Hopf operators on rearrangement-invariant
Banach function spaces}
Our main result is also flexible enough to treat the case of Wiener-Hopf operators on
many not necessarily separable rearrangement-invariant Banach function spaces.
\begin{theorem}
Let $X(\mathbb{R})$ be a rearrangement-invariant Banach function space with 
the lower Zippin index $p_X>0$ and let $X_2(\mathbb{R})$ be the closure of 
$L^2(\mathbb{R})\cap X(\mathbb{R})$ in the space $X(\mathbb{R})$.
 If $a\in\mathcal{M}_{X(\mathbb{R})}^0$, then the Wiener-Hopf operator 
 $W(a)\in\mathcal{B}(X_2(\mathbb{R}_+),X(\mathbb{R}_+))$
is maximally noncompact, that is,
\[
\|W(a)\|_{\mathcal{B}(X_2(\mathbb{R}_+),X(\mathbb{R}_+))}
=
\|W(a)\|_{\mathcal{B}(X_2(\mathbb{R}_+),X(\mathbb{R}_+)),\mathrm{e}}
=
\|W(a)\|_{\mathcal{B}(X_2(\mathbb{R}_+),X(\mathbb{R}_+)),\chi}.
\]
\end{theorem}
\begin{proof}
Following \cite[Ch.~3, Definition~1.2]{BS88}, for $1<p<\infty$, let 
$L^1(\mathbb{R})+L^p(\mathbb{R})$ be the collection of all measurable
functions $f:\mathbb{R}\to\mathbb{C}$ that are representable
as $f=g+h$ for some $g\in L^1(\mathbb{R})$ and $h\in L^p(\mathbb{R})$.
For each function $f\in L^1(\mathbb{R})+L^p(\mathbb{R})$, its norm is 
defined as 
\[
\|f\|_{L^1(\mathbb{R})+L^p(\mathbb{R})}
=
\inf\{\|g\|_{L^1(\mathbb{R})}+\|h\|_{L^p(\mathbb{R})}: \ f=g+h\},
\]
where the infimum is taken over all representations $f=g+h$ with
$g\in L^1(\mathbb{R})$ and $h\in L^p(\mathbb{R})$.

Since $p_X>0$, there exists $p\in(1,\infty)$ such that $p_X>1/p$. Then it 
follows from \cite[Lemma~5.2]{KS24} that $X(\mathbb{R})$ is continuously
embedded into $Z(\mathbb{R}):=L^1(\mathbb{R})+L^p(\mathbb{R})$.
Therefore $X(\mathbb{R_+})$ is contained in $Z(\mathbb{R}_+)$.
It is shown in the proof of \cite[Corollary~5.4]{KS24} that $S_0(\mathbb{R})$
is dense in $Z(\mathbb{R})$. Hence $S_0(\mathbb{R}_+)$ is dense
in $Z(\mathbb{R}_+)$ (cf. the proof of Lemma~\ref{le:density-simple}).
It remains to apply Theorem~\ref{th:main-result} to the translation-invariant
operator $A$ in \eqref{eq:Fourier-convolution-operator}, the Wiener-Hopf operator
$A_+=W(a)\in\mathcal{B}(X_2(\mathbb{R}_+),X(\mathbb{R}_+))$
and to the spaces $X_2(\mathbb{R})$, $X(\mathbb{R})$, and $Z(\mathbb{R})$ 
as above.
\end{proof}
Let us formulate a corollary of the above result, which does not follow from
Theorem~\ref{th:nice}.

For $1<p<\infty$, the Marcinkiewicz space $L^{p,\infty}(\mathbb{R})$
(the weak $L^p$ space) consists of all measurable functions 
$f:\mathbb{R}\to\mathbb{C}$ such that
\[
\|f\|_{L^{p,\infty}(\mathbb{R})}
:=
\sup_{0<t<\infty}t^{1/p-1}\int_0^t f^*(s)\,ds<\infty.
\]
It is well known that $L^{p,\infty}(\mathbb{R})$ is a rearrangement-invariant
Banach function space with respect to the above norm and both its Zippin indices
are equal to $1/p$ (see \cite[Ch.~4, Theorem~4.6]{BS88} and \cite[inequalities~(4.14)]{M85}).
Moreover, it follows from \cite[Ch.~II, Theorem 4.8 and Lemma 5.4]{KPS82}
or from \cite[Theorem 8.5.3]{PKJF13} and \cite[Ch.~1, Corollary~5.6]{BS88}
that $L^{p,\infty}(\mathbb{R})$ is nonseparable.
\begin{corollary}
Let $1<p<\infty$ and  let $L_2^{p,\infty}(\mathbb{R})$ be the closure of the set
$L^2(\mathbb{R})\cap L^{p,\infty}(\mathbb{R})$ in the Marcinkiewicz space 
$L^{p,\infty}(\mathbb{R})$.
 If $a\in\mathcal{M}_{L^{p,\infty}(\mathbb{R})}^0$, then the Wiener-Hopf operator 
 $W(a)\in\mathcal{B}(L_2^{p,\infty}(\mathbb{R}_+),L^{p,\infty}(\mathbb{R}_+))$
is maximally noncompact, that is,
\begin{align*}
\|W(a)\|_{\mathcal{B}(L_2^{p,\infty}(\mathbb{R}_+),L^{p,\infty}(\mathbb{R}_+))}
&=
\|W(a)\|_{\mathcal{B}(L_2^{p,\infty}(\mathbb{R}_+),L^{p,\infty}(\mathbb{R}_+)),\mathrm{e}}
\\
&=
\|W(a)\|_{\mathcal{B}(L_2^{p,\infty}(\mathbb{R}_+),L_2^{p,\infty}(\mathbb{R}_+)),\chi}.
\end{align*}
\end{corollary}
\section*{Declarations}
\subsection*{Acknowledgements}
We would like to thank the anonymous referees for useful comments and suggestions.

\subsection*{Funding}
This work is funded by national funds through the FCT -- Funda\c{c}\~ao para a 
Ci\^encia e a Tecnologia, I.P., under the scope of the projects UID/297/2025 
and UID/PRR/297/2020 
(Center for Mathematics and Applications -- NOVA Math).

\subsection*{Author contribution}
Both authors equally contributed to the manuscript.
\subsection*{Conflict of interest}
The authors declare no competing interests.
\subsection*{Data availability} 
The manuscript does not contain any associated data.
\bibliography{OKES25-Grudsky70}


\begin{thebibliography}{23}
\ifx \bisbn   \undefined \def \bisbn  #1{ISBN #1}\fi
\ifx \binits  \undefined \def \binits#1{#1}\fi
\ifx \bauthor  \undefined \def \bauthor#1{#1}\fi
\ifx \batitle  \undefined \def \batitle#1{#1}\fi
\ifx \bjtitle  \undefined \def \bjtitle#1{#1}\fi
\ifx \bvolume  \undefined \def \bvolume#1{\textbf{#1}}\fi
\ifx \byear  \undefined \def \byear#1{#1}\fi
\ifx \bissue  \undefined \def \bissue#1{#1}\fi
\ifx \bfpage  \undefined \def \bfpage#1{#1}\fi
\ifx \blpage  \undefined \def \blpage #1{#1}\fi
\ifx \burl  \undefined \def \burl#1{\textsf{#1}}\fi
\ifx \doiurl  \undefined \def \doiurl#1{\url{https://doi.org/#1}}\fi
\ifx \betal  \undefined \def \betal{\textit{et al.}}\fi
\ifx \binstitute  \undefined \def \binstitute#1{#1}\fi
\ifx \binstitutionaled  \undefined \def \binstitutionaled#1{#1}\fi
\ifx \bctitle  \undefined \def \bctitle#1{#1}\fi
\ifx \beditor  \undefined \def \beditor#1{#1}\fi
\ifx \bpublisher  \undefined \def \bpublisher#1{#1}\fi
\ifx \bbtitle  \undefined \def \bbtitle#1{#1}\fi
\ifx \bedition  \undefined \def \bedition#1{#1}\fi
\ifx \bseriesno  \undefined \def \bseriesno#1{#1}\fi
\ifx \blocation  \undefined \def \blocation#1{#1}\fi
\ifx \bsertitle  \undefined \def \bsertitle#1{#1}\fi
\ifx \bsnm \undefined \def \bsnm#1{#1}\fi
\ifx \bsuffix \undefined \def \bsuffix#1{#1}\fi
\ifx \bparticle \undefined \def \bparticle#1{#1}\fi
\ifx \barticle \undefined \def \barticle#1{#1}\fi
\bibcommenthead
\ifx \bconfdate \undefined \def \bconfdate #1{#1}\fi
\ifx \botherref \undefined \def \botherref #1{#1}\fi
\ifx \url \undefined \def \url#1{\textsf{#1}}\fi
\ifx \bchapter \undefined \def \bchapter#1{#1}\fi
\ifx \bbook \undefined \def \bbook#1{#1}\fi
\ifx \bcomment \undefined \def \bcomment#1{#1}\fi
\ifx \oauthor \undefined \def \oauthor#1{#1}\fi
\ifx \citeauthoryear \undefined \def \citeauthoryear#1{#1}\fi
\ifx \endbibitem  \undefined \def \endbibitem {}\fi
\ifx \bconflocation  \undefined \def \bconflocation#1{#1}\fi
\ifx \arxivurl  \undefined \def \arxivurl#1{\textsf{#1}}\fi
\csname PreBibitemsHook\endcsname

\bibitem[\protect\citeauthoryear{Lebow and Schechter}{1971}]{LS71}
\begin{barticle}
\bauthor{\bsnm{Lebow}, \binits{A.}},
\bauthor{\bsnm{Schechter}, \binits{M.}}:
\batitle{Semigroups of operators and measures of noncompactness}.
\bjtitle{J. Functional Analysis}
\bvolume{7},
\bfpage{1}--\blpage{26}
(\byear{1971})
\doiurl{10.1016/0022-1236(71)90041-3}
\end{barticle}
\endbibitem

\bibitem[\protect\citeauthoryear{Akhmerov et~al.}{1992}]{AKPRS92}
\begin{bbook}
\bauthor{\bsnm{Akhmerov}, \binits{R.R.}},
\bauthor{\bsnm{Kamenski\u{\i}}, \binits{M.I.}},
\bauthor{\bsnm{Potapov}, \binits{A.S.}},
\bauthor{\bsnm{Rodkina}, \binits{A.E.}},
\bauthor{\bsnm{Sadovski\u{\i}}, \binits{B.N.}}:
\bbtitle{Measures of Noncompactness and Condensing Operators}.
\bsertitle{Operator Theory: Advances and Applications},
vol. \bseriesno{55}.
\bpublisher{Birkh\"{a}user Verlag},
\blocation{Basel}
(\byear{1992}).
\doiurl{10.1007/978-3-0348-5727-7}
\end{bbook}
\endbibitem

\bibitem[\protect\citeauthoryear{Astala and Tylli}{1987}]{AT87}
\begin{barticle}
\bauthor{\bsnm{Astala}, \binits{K.}},
\bauthor{\bsnm{Tylli}, \binits{H.-O.}}:
\batitle{On the bounded compact approximation property and measures of
  noncompactness}.
\bjtitle{J. Funct. Anal.}
\bvolume{70}(\bissue{2}),
\bfpage{388}--\blpage{401}
(\byear{1987})
\doiurl{10.1016/0022-1236(87)90118-2}
\end{barticle}
\endbibitem

\bibitem[\protect\citeauthoryear{Krupnik}{2010}]{K10}
\begin{bchapter}
\bauthor{\bsnm{Krupnik}, \binits{N.}}:
\bctitle{Survey on the best constants in the theory of one-dimensional singular
  integral operators}.
In: \bbtitle{Topics in Operator Theory. {V}olume 1. {O}perators, Matrices and
  Analytic Functions}.
\bsertitle{Operator Theory: Advances and Applications},
vol. \bseriesno{202},
pp. \bfpage{365}--\blpage{393}.
\bpublisher{Birkh\"{a}user Verlag},
\blocation{Basel}
(\byear{2010})
\end{bchapter}
\endbibitem

\bibitem[\protect\citeauthoryear{Ayerbe~Toledano et~al.}{1997}]{ADl97}
\begin{bbook}
\bauthor{\bsnm{Ayerbe~Toledano}, \binits{J.M.}},
\bauthor{\bsnm{Dom\'{\i}nguez~Benavides}, \binits{T.}},
\bauthor{\bsnm{L\'{o}pez~Acedo}, \binits{G.}}:
\bbtitle{Measures of Noncompactness in Metric Fixed Point Theory}.
\bsertitle{Operator Theory: Advances and Applications},
vol. \bseriesno{99}.
\bpublisher{Birkh\"{a}user Verlag},
\blocation{Basel}
(\byear{1997}).
\doiurl{10.1007/978-3-0348-8920-9}
\end{bbook}
\endbibitem

\bibitem[\protect\citeauthoryear{Bana\'{s} and Mursaleen}{2014}]{BM14}
\begin{bbook}
\bauthor{\bsnm{Bana\'{s}}, \binits{J.}},
\bauthor{\bsnm{Mursaleen}, \binits{M.}}:
\bbtitle{Sequence Spaces and Measures of Noncompactness with Applications to
  Differential and Integral Equations}.
\bpublisher{Springer},
\blocation{New Delhi}
(\byear{2014}).
\doiurl{10.1007/978-81-322-1886-9}
\end{bbook}
\endbibitem

\bibitem[\protect\citeauthoryear{Luxemburg}{1955}]{L55}
\begin{bbook}
\bauthor{\bsnm{Luxemburg}, \binits{W.A.J.}}:
\bbtitle{Banach Function Spaces}.
\bpublisher{Technische Hogeschool te Delft},
\blocation{Delft}
(\byear{1955})
\end{bbook}
\endbibitem

\bibitem[\protect\citeauthoryear{Bennett and Sharpley}{1988}]{BS88}
\begin{bbook}
\bauthor{\bsnm{Bennett}, \binits{C.}},
\bauthor{\bsnm{Sharpley}, \binits{R.}}:
\bbtitle{Interpolation of Operators}.
\bsertitle{Pure and Applied Mathematics},
vol. \bseriesno{129}.
\bpublisher{Academic Press},
\blocation{Boston, MA}
(\byear{1988})
\end{bbook}
\endbibitem

\bibitem[\protect\citeauthoryear{Pick et~al.}{2013}]{PKJF13}
\begin{bbook}
\bauthor{\bsnm{Pick}, \binits{L.}},
\bauthor{\bsnm{Kufner}, \binits{A.}},
\bauthor{\bsnm{John}, \binits{O.}},
\bauthor{\bsnm{Fu\v{c}\'{\i}k}, \binits{S.}}:
\bbtitle{Function Spaces. {V}ol. 1},
\bedition{extended} edn.
\bsertitle{De Gruyter Series in Nonlinear Analysis and Applications},
vol. \bseriesno{14}.
\bpublisher{Walter de Gruyter \& Co.},
\blocation{Berlin}
(\byear{2013})
\end{bbook}
\endbibitem

\bibitem[\protect\citeauthoryear{Maz'ya and Shaposhnikova}{2009}]{MS09}
\begin{bbook}
\bauthor{\bsnm{Maz'ya}, \binits{V.G.}},
\bauthor{\bsnm{Shaposhnikova}, \binits{T.O.}}:
\bbtitle{Theory of {S}obolev Multipliers. With Applications to Differential and
  Integral Operators}.
\bsertitle{Grundlehren der mathematischen Wissenschaften},
vol. \bseriesno{337}.
\bpublisher{Springer},
\blocation{Berlin}
(\byear{2009}).
\doiurl{10.1007/978-3-540-69492-2}
\end{bbook}
\endbibitem

\bibitem[\protect\citeauthoryear{Edmunds and Shargorodsky}{2005}]{ES05}
\begin{barticle}
\bauthor{\bsnm{Edmunds}, \binits{D.E.}},
\bauthor{\bsnm{Shargorodsky}, \binits{E.}}:
\batitle{The inner variation of an operator and the essential norms of
  pointwise multipliers in function spaces}.
\bjtitle{Houston J. Math.}
\bvolume{31}(\bissue{3}),
\bfpage{841}--\blpage{855}
(\byear{2005})
\end{barticle}
\endbibitem

\bibitem[\protect\citeauthoryear{Karlovych and Shargorodsky}{2024}]{KS24}
\begin{barticle}
\bauthor{\bsnm{Karlovych}, \binits{O.}},
\bauthor{\bsnm{Shargorodsky}, \binits{E.}}:
\batitle{Discrete {R}iesz transforms on rearrangement-invariant {B}anach
  sequence spaces and maximally noncompact operators}.
\bjtitle{Pure Appl. Funct. Anal.}
\bvolume{9}(\bissue{1}),
\bfpage{195}--\blpage{210}
(\byear{2024})
\end{barticle}
\endbibitem

\bibitem[\protect\citeauthoryear{Karlovych and Shargorodsky}{2022}]{KS22}
\begin{barticle}
\bauthor{\bsnm{Karlovych}, \binits{O.}},
\bauthor{\bsnm{Shargorodsky}, \binits{E.}}:
\batitle{On the essential norms of singular integral operators with constant
  coefficients and of the backward shift}.
\bjtitle{Proc. Amer. Math. Soc. Ser. B}
\bvolume{9},
\bfpage{60}--\blpage{70}
(\byear{2022})
\doiurl{10.1090/bproc/118}
\end{barticle}
\endbibitem

\bibitem[\protect\citeauthoryear{Edmunds and Lang}{}]{EL25-1}
\begin{botherref}
\oauthor{\bsnm{Edmunds}, \binits{D.E.}},
\oauthor{\bsnm{Lang}, \binits{J.}}:
Remarks on the {L}aplace transfrom.
Rev. Mat. Complut.
\doiurl{10.1007/s13163-025-00542-8}
\end{botherref}
\endbibitem

\bibitem[\protect\citeauthoryear{Edmunds and Lang}{2025}]{EL25-2}
\begin{barticle}
\bauthor{\bsnm{Edmunds}, \binits{D.E.}},
\bauthor{\bsnm{Lang}, \binits{J.}}:
\batitle{Notes on non-compact maps and the importance of {B}ernstein numbers}.
\bjtitle{Adv. Oper. Theory}
\bvolume{10}(\bissue{4}),
\bfpage{86}--\blpage{30}
(\byear{2025})
\doiurl{10.1007/s43036-025-00456-8}
\end{barticle}
\endbibitem

\bibitem[\protect\citeauthoryear{B\"{o}ttcher et~al.}{1997}]{BGS97}
\begin{barticle}
\bauthor{\bsnm{B\"{o}ttcher}, \binits{A.}},
\bauthor{\bsnm{Grudsky}, \binits{S.M.}},
\bauthor{\bsnm{Silbermann}, \binits{B.}}:
\batitle{Norms of inverses, spectra, and pseudospectra of large truncated
  {W}iener-{H}opf operators and {T}oeplitz matrices}.
\bjtitle{New York J. Math.}
\bvolume{3},
\bfpage{1}--\blpage{31}
(\byear{1997})
\end{barticle}
\endbibitem

\bibitem[\protect\citeauthoryear{B\"{o}ttcher and Silbermann}{2006}]{BS06}
\begin{bbook}
\bauthor{\bsnm{B\"{o}ttcher}, \binits{A.}},
\bauthor{\bsnm{Silbermann}, \binits{B.}}:
\bbtitle{Analysis of {T}oeplitz Operators},
\bedition{2}nd edn.
\bsertitle{Springer Monographs in Mathematics}.
\bpublisher{Springer},
\blocation{Berlin}
(\byear{2006}).
\doiurl{10.1007/3-540-32436-4}
\end{bbook}
\endbibitem

\bibitem[\protect\citeauthoryear{Karlovych and Valente}{2025}]{KV25}
\begin{bchapter}
\bauthor{\bsnm{Karlovych}, \binits{O.}},
\bauthor{\bsnm{Valente}, \binits{M.}}:
\bctitle{On the operator and essential norms of {F}ourier convolution operators
  and {W}iener-{H}opf operators with the same symbol}.
In: \bbtitle{Operator Theory, Related Fields, and Applications}.
\bsertitle{Operator Theory: Advances and Applications},
vol. \bseriesno{307},
pp. \bfpage{361}--\blpage{377}.
\bpublisher{Birkh\"{a}user/Springer},
\blocation{Cham}
(\byear{2025})
\end{bchapter}
\endbibitem

\bibitem[\protect\citeauthoryear{Cruz-Uribe and Fiorenza}{2013}]{CF13}
\begin{bbook}
\bauthor{\bsnm{Cruz-Uribe}, \binits{D.V.}},
\bauthor{\bsnm{Fiorenza}, \binits{A.}}:
\bbtitle{Variable {L}ebesgue Spaces. Foundations and Harmonic Analysis}.
\bsertitle{Applied and Numerical Harmonic Analysis}.
\bpublisher{Birkh\"{a}user/Springer},
\blocation{Heidelberg}
(\byear{2013}).
\doiurl{10.1007/978-3-0348-0548-3}
\end{bbook}
\endbibitem

\bibitem[\protect\citeauthoryear{Lorist and Nieraeth}{2024}]{LN24}
\begin{barticle}
\bauthor{\bsnm{Lorist}, \binits{E.}},
\bauthor{\bsnm{Nieraeth}, \binits{Z.}}:
\batitle{Banach function spaces done right}.
\bjtitle{Indag. Math. (N.S.)}
\bvolume{35}(\bissue{2}),
\bfpage{247}--\blpage{268}
(\byear{2024})
\doiurl{10.1016/j.indag.2023.11.004}
\end{barticle}
\endbibitem

\bibitem[\protect\citeauthoryear{Karlovich and Shargorodsky}{2019}]{KS19}
\begin{barticle}
\bauthor{\bsnm{Karlovich}, \binits{A.}},
\bauthor{\bsnm{Shargorodsky}, \binits{E.}}:
\batitle{When does the norm of a {F}ourier multiplier dominate its {$L^\infty$}
  norm?}
\bjtitle{Proc. Lond. Math. Soc. (3)}
\bvolume{118}(\bissue{4}),
\bfpage{901}--\blpage{941}
(\byear{2019})
\doiurl{10.1112/plms.12206}
\end{barticle}
\endbibitem

\bibitem[\protect\citeauthoryear{Maligranda}{1985}]{M85}
\begin{barticle}
\bauthor{\bsnm{Maligranda}, \binits{L.}}:
\batitle{Indices and interpolation}.
\bjtitle{Dissertationes Math. (Rozprawy Mat.)}
\bvolume{234},
\bfpage{49}
(\byear{1985})
\end{barticle}
\endbibitem

\bibitem[\protect\citeauthoryear{Kre\u{\i}n et~al.}{1982}]{KPS82}
\begin{bbook}
\bauthor{\bsnm{Kre\u{\i}n}, \binits{S.G.}},
\bauthor{\bsnm{Petun\={\i}n}, \binits{Y.I.}},
\bauthor{\bsnm{Sem\"{e}nov}, \binits{E.M.}}:
\bbtitle{Interpolation of Linear Operators}.
\bsertitle{Translations of Mathematical Monographs},
vol. \bseriesno{54}.
\bpublisher{American Mathematical Society},
\blocation{Providence, RI}
(\byear{1982})
\end{bbook}
\endbibitem

\end{thebibliography}
\end{document}